\documentclass[10pt]{article}

\usepackage{amsmath,amssymb,amsthm,amsfonts,amscd,flafter,epsf, epsfig,graphicx}

\title{A note on genus one fibered knots in lens spaces}
\author{John A. Baldwin}
\date{}

\newcommand{\B}{\ensuremath{{\mathfrak b}}}
\newcommand\by{\delta}
\newcommand\w{x y^2 x y^2}

\newtheorem{thm}{Theorem}[section]

\newtheorem{cor}[thm]{Corollary}

\begin{document}
\maketitle
\begin{abstract}  
Supose that $Y$ is a lens space with $|H_1(Y;\mathbb{Z})|$ prime, and $Y$ does not contain a genus one fibered knot. We show that $Y$ contains a knot whose exterior is a once-punctured torus bundle if and only if $Y$ is the result of $p/q$-surgery on the trefoil. This partially answers a question posed by Ken Baker in a paper in which he gives a complete classification of genus one fibered knots contained in lens spaces. Combining Baker's classification with Moser's characterization of lens space surgeries on the trefoil, we generate an infinite family of lens spaces which do not contain any knot whose exterior is a once-punctured torus bundle.
 \end{abstract} 

\section{Introduction}
When $Y$ is a three-manifold we say that a knot $K \subset Y$ is a \emph{genus one fibered knot} (GOF-knot) if

\begin{enumerate}
\item{$Y - \nu(K)$ is a once-punctured torus bundle (OPTB),}
\item{the monodromy restricts to the identity on the boundary of the fiber, and}
\item{$K$ is ambient isotopic in $Y$ to the boundary of a fiber}
\end{enumerate}

\noindent Given a three-manifold, one might ask how many distinct GOF-knots it contains. For instance, it is well known that the only GOF-knots in $S^3$ are the right- and left-handed trefoil knots and the figure eight knot. Morimoto partially answers this question for lens spaces and, in the process, shows that $L(19,2)$, $L(19,4)$, and $L(19,7)$ do not contain any GOF-knots \cite{Mor}. In the same paper, Morimoto asks whether there is an upper limit to the number of GOF-knots that may be contained in a single lens space. 

In early 2006, Ken Baker answered this question in the affirmative, considerably extending Morimoto's results \cite{Bak}. Baker shows that the lens spaces $L(m,1)$, where $m>0$, contain exactly two GOF-knots, except for $L(4,1)$ which contains three. He proves that all other lens spaces contain at most one GOF-knot. More generally, Baker gives a criterion for determining exactly how many GOF-knots a given lens space contains.

Using this criterion, it is possible to find a plethora of lens spaces with no GOF-knots. For instance, one can verify that the examples first discovered by Morimoto, $L(19,2)$, $L(19,4)$, and $L(19,7)$, do not contain any GOF-knots. In light of the fact that there are lens spaces with no GOF-knots, Baker asks whether every lens space must contain a knot whose exterior is a OPTB? For instance, although $L(19,7)$ contains no GOF-knots, it is the result of $-19/3$-surgery on the left-handed trefoil $T(3,2)$, and therefore contains a knot whose exterior is a OPTB, namely the core of the surgered torus. We aim to address this question. Our main result is the following:

\begin{thm}
\label{thm:MainTheorem}
Suppose that $Y$ is a three-manifold which does not contain a GOF-knot. Suppose further that $|H_1(Y, \mathbb{Z})|$ is prime. Then $Y$ contains a knot whose exterior is a OPTB if and only if $Y$ can be obtained as $p/q$-surgery on the trefoil, or on the figure eight, where $p>1$. 
\end{thm}

\noindent Note that if $Y$ is a lens space (or, more generally, a Heegaard Floer homology $L$-space), then it is well known that $Y$ cannot be obtained by surgery on the figure eight. See, for example, \cite{OSz}.

Baker's approach relies on the fact that a three-manifold $Y$ which contains a GOF-knot $K$ is the branched double cover of $S^3$, branched over a closed 3-braid. In this formulation, $K$ is the lift of the braid axis of this 3-braid. In the case of lens spaces, it is known that $L(\alpha, \beta)$ is the double-cover of $S^3$ branched over a link $L$ if and only if $L$ is equivalent to the 2-bridge link $\B(\alpha, \beta)$. Along these lines, Baker proves the following theorem:

\begin{thm}
\label{thm:BakerTheorem}
The lens space $L(\alpha, \beta)$ contains exactly $N$ GOF-knots if and only if the 2-bridge link $\B(\alpha,\beta)$ admits exactly $N$ equivalence classes of 3-braid representatives. No lens space contains four GOF-knots. Furthermore if $L$ is an unoriented 2-bridge link, then
\begin{enumerate}
\item{$L$ admits three equivalence classes of 3-braid representatives only if $L$ is equivalent to $\B(4,1)$.}
\item{$L$ admits two equivalence classes of 3-braid representatives only if $L$ is equivalent to $\B(\alpha,1)$, and $\alpha \neq 0$.}
\item{$L$ admits exactly one 3-braid representatives if $L$ is equivalent to either $\B(0,1)$, or $\B(\alpha, \beta)$ where $0<\beta<\alpha$ and either
\begin{itemize}
\item{$\alpha = 2pq+p+q$ and $\beta = 2q+1$ for some integers $p,q>1$, or}
\item{$\alpha = 2pq+p+q+1$ and $\beta = 2q+1$ for some integers $p,q>0$.}
\end{itemize}
}
\item{$L$ admits no 3-braid representatives otherwise.}
\end{enumerate}
\end{thm}

\subsection{Organization}
In section \ref{sec:Proof}, we recount our previous result on the classification of genus one, one boundary component open books in \cite{Bal} and use this to prove Theorem \ref{thm:MainTheorem}. In section \ref{sec:Examples}, we combine Moser's results on lens space surgeries \cite{Mos} with Theorems \ref{thm:MainTheorem} and \ref{thm:BakerTheorem} to produce an infinite family of lens spaces which do not contain any knot whose exterior is a OPTB.

\section{Proof of Theorem \ref{thm:MainTheorem}}
\label{sec:Proof}
The mapping class group of the once-punctured torus $\Sigma$ is generated by right-handed Dehn twists about dual non-separating curves, $x$ and $y$. In an abuse of notation we denote, by $\gamma$, the right-handed Dehn twist around the curve $\gamma \subset \Sigma.$ The left-handed Dehn twist around $\gamma$ is then denoted by $\gamma^{-1}$. Thus, given an open-book decomposition ($\Sigma$, $\phi$), we can express $\phi$ as a product of Dehn twists, $x^{a_1} y^{b_1}x^{a_2} y^{b_2}...x^{a_n} y^{b_n}$, with $a_{i}, b_{j} \in \mathbb{Z}$ (In our notation, composition is on the left). Using relations in the mapping class group of $\Sigma$, we can separate the genus one, one boundary component open books into six categories according to the following \cite{Bal}:

\begin {thm}
\label{thm:OBClassification}
Let $\delta$ be a curve in $\Sigma$ parallel to the boundary. Then any genus one, one boundary component open book can be written as $(\Sigma, \phi)$, where $\phi$ is one of the following types:
\begin{description}
\item[A.] $\by^d *x^{a_1}y^{-1}...x^{a_n}y^{-1}$, where the $a_i \geq 0$, some $a_j \neq 0$.
\item[B.] $\by^d * \w*x^{a_1}y^{-1}...x^{a_n}y^{-1}$, where the $a_i \geq 0$, some $a_j \neq 0$.
\item[C.] $\by^d *y^m$, for $m \in \mathbb{Z}$.
\item[D.] $\by^d *\w*y^m$, for $m \in \mathbb{Z}$.
\item[E.] $\by^d *x^m y^{-1}$, where $m \in \{-1,-2,-3\}$
\item[F.] $\by^d *\w*x^m y^{-1}$, where $m \in \{-1,-2,-3\}$
\end{description}
\end{thm}

\begin{figure}[!htbp]
\begin{center}
\includegraphics[height=8cm]{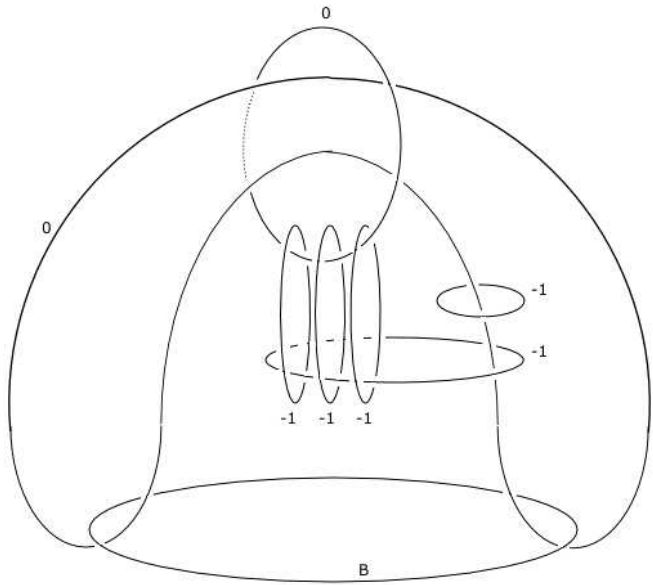}
\caption{\quad A surgery diagram for the open book given by the monodromy $\phi = x y^3 x$. The curve $B$ is the binding.}
\label{fig:KirbyDiagram}
\end{center}
\end{figure}

In the surgery diagram for these open books (see Figure \ref{fig:KirbyDiagram} for an example), the binding $B$ (the fibered knot corresponding to this open book) does not algebraically link the other surgery curves. 

\noindent Hence, if $M$ has an open book decomposition and $M_{p/q}(B)$ denotes the result of $p/q$-surgery on $B$, then for $p>1$ 

$$H_1(M_{p/q}(B), \mathbb{Z}) = H_1(M, \mathbb{Z}) \oplus \mathbb{Z} /p\mathbb{Z}.$$ 

Suppose now that $Y$ is a three-manifold which does not contain a GOF-knot, but which does contain a knot $K$ whose exterior is a OPTB. Then $Y$ must be some $p/q$-Dehn-filling of a OPTB with boundary-fixing monodromy (if necessary, we can introduce twists around a curve parallel to the boundary of the fiber surface). Therefore, we can express $Y$ as $p/q$-surgery on the binding of some genus one, one boundary component open book, where $p>1$ (If $p = 1$, then $Y$ has a genus one, one boundary component open book decomposition, and therefore contains a fibered knot). We write $Y = M_{p/q}(B)$. Now suppose that $|H_1(Y, \mathbb{Z})|$ is prime. Then it follows from the discussion in the previous paragraph that $H_1(M, \mathbb{Z})=0$. This allows us to narrow down the possibilities for $M$.

For the monodromies listed in Theorem \ref{thm:OBClassification}, it is relatively easy to calculate the first homology of their corresponding open books - we merely have to compute the determinant of the linking matrices for the associated surgery diagrams. The Appendix of \cite{Bal} contains a computation of $|H_1|$ for the open books in Theorem \ref{thm:OBClassification}.B. For these open books, it is shown that $|H_1|$ is positive and strictly increasing in the $a_i$. One can show, more or less identically, that the same is true for the open books in Theorem \ref{thm:OBClassification}.A. The computation of $|H_1|$ for the open books in Theorem \ref{thm:OBClassification}.C-F is very simple as the surgery diagram is uncomplicated. We summarize all of this below:

\begin{description}
\item[A.] $|H_1| = 1$ for $n=1$ and $a_1 =1$; $|H_1|>1$ otherwise.
\item[B.] $|H_1| \geq 4$.
\item[C.] $|H_1| = \infty$.
\item[D.] $|H_1| = 4$.
\item[E.] $|H_1| = 1$ for $m=-1$; $|H_1|>1$ otherwise.
\item[F.] $|H_1| = 1$ for $m=-3$; $|H_1|>1$ otherwise.
\end{description}

Therefore, the open books which are candidates for $M$ in the above discussion are given by the monodromies 

\begin{enumerate}
\item{$\by^d * x y^{-1}$}
\item{$\by^d * x^{-1} y^{-1}$}
\item{$\by^d * \w * x^{-3} y^{-1} =\by^d * x y$}
\end{enumerate}

\noindent $p/q$-surgery on the bindings of these open books yields $p/(q-dp)$-surgery on (respectively) 
\begin{enumerate}
\item{the binding of the open book $(\Sigma, xy^{-1})$, i.e. the figure eight.}
\item{the binding of the open book $(\Sigma, x^{-1}y^{-1})$, i.e. the left-handed trefoil.}
\item{the binding of the open book $(\Sigma, xy)$, i.e. the right-handed trefoil.}
\end{enumerate}

This proves Theorem \ref{thm:MainTheorem}. As mentioned before, if $Y$ is a lens space, $Y$ cannot be obtained as surgery on the figure eight. \qed

\section{Examples}
\label{sec:Examples}
In order to find a lens space which does not contain any knot whose exterior is a OPTB, we have shown that it is sufficient to find a lens space $L(m,n)$ such that $m$ is prime, $L(m,n)$ does not contain a GOF-knot, and $L(m,n)$ is not the result of surgery on the trefoil. Baker gives a classification of lens spaces with no GOF-knots, and Moser characterizes precisely those surgeries on the torus knot $T(r,s)$ which yield lens spaces. Moser proves the following theorem \cite{Mos}:

\begin{thm}
\label{thm:MoserTheorem}
$p/q$-surgery on the torus knot $T(r,s)$ for $0<s<r$ yields a lens space if and only if $|rsq+p| = 1$. The lens space produced by such a surgery is $L(|p|,qs^2)$.
\end{thm}

\noindent Hence, for the left-handed trefoil, which is the torus knot $T(3,2)$, this condition becomes $|6q+p|=1$, and the lens space obtained is $L(|p|, 4q)$. This gives us a means to find many lens spaces with no knot whose exterior is a OPTB. The following is a corollary of Theorems \ref{thm:BakerTheorem}, \ref{thm:MainTheorem}, and \ref{thm:MoserTheorem} and provides one infinite family of such lens spaces.

\begin{cor}
\label{cor:Examples}
The lens spaces $L(m,2)$ contain no knot whose exterior is a OPTB when $m$ is prime and $m \geq 11$.
\end{cor}

\begin{proof}
The lens spaces $L(\alpha,\beta)$ with $0<\beta<\alpha$ which are homeomorphic to $L(m,2)$ are $L(m,2)$, $L(m,m-2)$, $L(m, (m+1)/2)$, and $L(m, (m-1)/2)$. Suppose that $p/q$-surgery on $T(3,2)$ yields $L(m,2)$. Then $m = |p|$ and, by Moser, $|6q+p|=1$, which implies that $12q \equiv \pm 2\ mod\ m$. Also, by Moser, $L(m,2)$ is homeomorphic to $L(m, 4q)$. Hence, $4q \equiv 2, -2, (m+1)/2,$ or $(m-1)/2\ mod\ m$. Multiplying the first equation by two and the second by six, we get

\begin{eqnarray*}
&&24q \equiv \pm 4\ mod\ m \\
&&24q \equiv 12, -12, 3, or-3\ mod\ m \\
\end{eqnarray*}

\noindent But these equations are inconsistent if $m \geq 11$ is prime. To complete the proof, we need to show that $L(m,2)$ contains no GOF-knots. According to Theorem \ref{thm:BakerTheorem}, it is sufficient to show that if $2k+1 = 2,m-2,(m+1)/2,$ or $(m-1)/2$ for some $k >0$, then $2kl+k+l \neq m$ or $m-1$ for any $l >0$.  

\begin{itemize}
\item{If $2k+1 = m-2$, then $2kl+k+l = l(m-2)+(m-3)/2$ is never equal to $m$ or $m-1$ when $m>7$.}
\item{If $2k+1 = (m+1)/2$, then $2kl+k+l = l(m+1)/2+(m-1)/4$ is never equal to $m$ or $m-1$ when $m>5$.}
\item{If $2k+1 = (m-1)/2$, then $2kl+k+l = l(m-1)/2+(m-3)/4$ is never equal to $m$ or $m-1$ when $m>3$.}
\end{itemize}

\noindent Altogether, we've shown that $L(m,2)$ for $m \geq 11$ and prime does not contain any GOF-knots and cannot be obtained by surgery on the trefoil. Thus, by Theorem \ref{thm:MainTheorem}, $L(m,2)$ contains no knot whose exterior is a OPTB.
\end{proof}


\begin{thebibliography}{}

\bibitem[Bak]{Bak}
K. Baker. \textit{Counting genus one fibered knots in lens spaces}. preprint, 2006.

\bibitem[Bal]{Bal}
J. A. Baldwin. \textit{Tight contact structures and genus one fibered knots}. preprint, 2006.

\bibitem[Mor]{Mor}
K. Morimoto. \textit{Genus one fibered knots in lens spaces}. J.\ Math.\ Soc.\ Japan\ {\bf 41}(1) (1989), 81-96.

\bibitem[Mos]{Mos}
L. Moser. \textit{Elementary surgery along a torus knot}. Pac.\ J.\ Math.\ {\bf 41} (1971), 737-745.

\bibitem[OSz]{OSz}
P. Ozsv{\'a}th and Z. Szab{\'o}. \textit{On knot Floer homolgy and lens space surgeries}. preprint, 2004.

\end{thebibliography}
\end{document}